\documentclass[10pt]{article}
\usepackage{cite}
\usepackage{mathrsfs}
\usepackage{amsfonts}
\usepackage{amsmath}
\usepackage{amsfonts,amssymb}
\usepackage{dsfont}
\usepackage{curves}
\usepackage{mathrsfs}
\usepackage{pifont}
\usepackage{amssymb}
\allowdisplaybreaks

\numberwithin{equation}{section}

\date{}

\def\BigRoman{\uppercase\expandafter{\romannumeral\number\count 255 }}
\def\Romannumeral{\afterassignment\BigRoman\count255=}

\begin{document}
\title{Spanning subgraphs and spectral radius in graphs
}
\author{\small  Sizhong Zhou\footnote{Corresponding
author. E-mail address: zsz\_cumt@163.com (S. Zhou)}\\
\small School of Science, Jiangsu University of Science and Technology,\\
\small Zhenjiang, Jiangsu 212100, China\\
}

\maketitle
\begin{abstract}
\noindent A spanning tree $T$ of a connected graph $G$ is a subgraph of $G$ that is a tree covers all vertices of $G$. The leaf distance of $T$ is
defined as the minimum of distances between any two leaves of $T$. A fractional matching of a graph $G$ is a function $h$ assigning every edge a
real number in $[0,1]$ so that $\sum\limits_{e\in E_G(v)}{h(e)}\leq1$ for any $v\in V(G)$, where $E_G(v)$ denotes the set of edges incident with $v$
in $G$. A fractional matching of $G$ is called a fractional perfect matching if $\sum\limits_{e\in E_G(v)}{h(e)}=1$ for any $v\in V(G)$. A graph $G$
with at least $2k+2$ vertices is said to be fractional $k$-extendable if every $k$-matching $M$ in $G$ is included in a fractional perfect matching
$h$ of $G$ such that $h(e)=1$ for any $e\in M$. This paper considers a lower bound on the spectral radius of $G$ to guarantee that $G$ has a spanning
tree with leaf distance at least $d$. At the same time, we obtain a lower bound on the spectral radius of $G$ to ensure that $G$ is fractional
$k$-extendable.
\\
\begin{flushleft}

(2020) Mathematics Subject Classification: 05C50, 05C05, 05C70
\end{flushleft}
\end{abstract}

\section{Introduction}

In this paper, all graphs considered are assumed to simple and undirected. Let $G=(V(G),E(G))$ denote a graph with vertex set $V(G)$ and edge
set $E(G)$. The order and size of $G$ is denoted by $n$ and $e(G)$, respectively. That is, $n=|V(G)|$ and $e(G)=|E(G)|$. For $v\in V(G)$, the
set of vertices adjacent to $v$ in $G$ is called the neighborhood of $v$ and denoted by $N_G(v)$. We denote by $d_G(v)=|N_G(v)|$ the degree of
$v$ in $G$, and by $\delta(G)$ (or $\delta$ for short) the minimum degree of $G$. Let $\alpha(G)$ and $i(G)$ denote the independence number and
the number of isolated vertices in $G$, respectively. For any $S\subseteq V(G)$, let $G[S]$ denote the subgraph of $G$ induced by $S$, and write
$G-S=G[V(G)\setminus S]$. The complete graph of order $n$ is denoted by $K_n$. Let $c$ be a real number. Recall that $\lceil c\rceil$ is the
smallest integer satisfying $\lceil c\rceil\geq c$.

Let $G_1$ and $G_2$ be two vertex-disjoint graphs. The union of $G_1$ and $G_2$ is denoted by $G_1\cup G_2$, which is the graph with vertex set
$V(G_1)\cup V(G_2)$ and edge set $E(G_1)\cup E(G_2)$. The join $G_1\vee G_2$ is obtained from $G_1\cup G_2$ by adding all the edges joining a
vertex of $G_1$ to a vertex of $G_2$.

Given a graph $G$ with vertex set $V(G)=\{v_1,v_2,\ldots,v_n\}$, let $A(G)$ denote the adjacency matrix of $G$. The $(i,j)$-entry of $A(G)$ is 1
if $v_iv_j\in E(G)$, and 0 otherwise. The eigenvalues of $A(G)$ are called the eigenvalues of $G$. It is obvious that $A(G)$ is a real symmetric
nonnegative matrix. Consequently, its eigenvalues are real, which can be arranged in non-increasing order as
$\lambda_1(G)\geq\lambda_2(G)\geq\cdots\geq\lambda_n(G)$. Notice that the spectral radius of $G$, denoted by $\rho(G)$, is equal to $\lambda_1(G)$.

A spanning tree $T$ of a connected graph $G$ is a subgraph of $G$ that is a tree covers all vertices of $G$. For $v\in V(T)$, the vertex $v$ is
called a leaf of $T$ if $d_T(v)=1$. The leaf degree of a vertex $v\in V(T)$ is defined as the number of leaves adjacent to $v$ in $T$. The leaf
degree of $T$ is the maximum leaf degree among all the vertices of $T$. The leaf distance of $T$ is defined as the minimum of distances between
any two leaves of $T$. In fact, a tree with leaf degree 1 has leaf distance at least 3.

Kaneko \cite{Ks} presented some sufficient conditions for a connected graph to have a spanning tree with leaf distance at least $d=3$ and
conjectured that similar conditions suffice for larger $d$. Later, Kaneko, Kano and Suzuki \cite{KKS} claimed that Kaneko's conjecture is true
for $d=4$. For $d\geq4$, Erbes, Molla, Mousley and Santana \cite{EMMS} showed that a stronger form of Kaneko's conjecture holds for all $n$-vertex
connected graphs with $\alpha(G)\leq5$, and proved Kaneko's conjecture for $d\geq\frac{n}{3}$. Zhou, Sun and Liu \cite{ZSL3} provided two spectral
conditions for a connected graph to contain a spanning tree with leaf distance at least $d=3$. Chen, Lv, Li and Xu \cite{CLLX} investigated the
existence of spanning trees with leaf distance at least $d=4$ in connected graphs and obtained three new results. More results on spanning trees
can be found in \cite{GS,Kyaw,MM,ZW,ZZL,Wc}.

A set $M\subseteq E(G)$ is a matching if no two edges share a vertex. A matching of size $k$ is called a $k$-matching. A matching $M$ is called a
perfect matching (or 1-factor) if it covers all the vertices of $G$. Let $k\geq0$ be an integer. Then a graph $G$ with at least $2k+2$ vertices is
said to be $k$-extendable if every $k$-matching in $G$ can be extended to a perfect matching in $G$. A fractional matching of a graph $G$ is a
function $h$ assigning every edge a real number in $[0,1]$ so that $\sum\limits_{e\in E_G(v)}{h(e)}\leq1$ for any $v\in V(G)$, where $E_G(v)$
denotes the set of edges incident with $v$ in $G$. A fractional matching of $G$ is called a fractional perfect matching if
$\sum\limits_{e\in E_G(v)}{h(e)}=1$ for any $v\in V(G)$. Then a graph $G$ with at least $2k+2$ vertices is said to be fractional $k$-extendable if
every $k$-matching $M$ in $G$ is included in a fractional perfect matching $h$ of $G$ such that $h(e)=1$ for any $e\in M$.

The perfect matching and matching extendability attracted much attention. Tutte \cite{T} provided a characterization for a graph to contain a perfect
matching. Enomoto \cite{E} established a connection between toughness and a perfect matching in a graph. Niessen \cite{N} presented a neighborhood
union condition for a graph to have a perfect matching. O \cite{Os} obtained a spectral radius condition to guarantee that a connected graph has a
perfect matching. Zhang and Lin \cite{ZLp} got a distance spectral condition to guarantee the existence of a perfect matching in a connected graph.
Plummer \cite{P1} first introduced the concept of $k$-extendable graph and obtained some results on $k$-extendable graphs. Ananchuen and Caccetta
\cite{AC}, Lou and Yu \cite{LY}, Cioaba, Koolen and Li \cite{CKL}, Robertshaw and Woodall \cite{RW} investigated the existence of $k$-extendable
graphs. The fractional perfect matching and fractional matching extendability also attracted much attention. Lov\'asz and Plummer \cite{LP} showed
a characterization for the existence of fractional perfect matchings in graphs. Liu and Zhang \cite{LZ} claimed a toughness condition for a graph
to contain a fractional perfect matching in a graph. Ma and Liu \cite{ML} provided a characterization of fractional $k$-extendable graphs. Zhu and
Liu \cite{ZLs} established a relationship between binding numbers and fractional $k$-extendable graphs. Much effort has been devoted to finding
sufficient conditions for the existence of spanning subgraphs (see \cite{Zs,ZSL1,Zr,ZZL2,WZhi,Wp,M,ZWa,ZZS,GWC,Zt,ZZL1}).

Motivated by \cite{Os,EMMS,ML} directly, we are to establish a spectral radius condition for the existence of a spanning tree with leaf distance at
least $d$ in a connected graph, and propose a lower bound on the spectral radius of a connected graph $G$ to guarantee that $G$ is a fractional
$k$-extendable graph. Our main results are shown in the following.

\medskip

\noindent{\textbf{Theorem 1.1.}} Let $G$ be a connected graph of order $n$ with $\alpha(G)\leq5$, and let $d$ be an integer with $16\leq d^{2}\leq n$.
If
$$
\rho(G)\geq\rho(K_{\lceil\frac{d}{2}\rceil-1}\vee(K_{n-\lceil\frac{d}{2}\rceil}\cup K_1)),
$$
then $G$ has a spanning tree with leaf distance at least $d$, unless $G=K_{\lceil\frac{d}{2}\rceil-1}\vee(K_{n-\lceil\frac{d}{2}\rceil}\cup K_1)$.

\medskip

\noindent{\textbf{Theorem 1.2.}} Let $k\geq1$ be an integer, and let $G$ be a connected graph of order $n$ with minimum degree $\delta$ and
$n\geq\max\{2k+9,5\delta+1\}$. If
$$
\rho(G)\geq\max\{\rho(K_{2k}\vee(K_{n-2k-1}\cup K_1)), \rho(K_{\delta}\vee(K_{n-2\delta+2k-1}\cup(\delta-2k+1)K_1))\},
$$
then $G$ is fractional $k$-extendable, unless $G\in\{K_{2k}\vee(K_{n-2k-1}\cup K_1),K_{\delta}\vee(K_{n-2\delta+2k-1}\cup(\delta-2k+1)K_1)\}$.

\medskip

The proofs of Theorems 1.1 and 1.2 will be provided in Sections 3 and 4, respectively.

\section{Preliminary lemmas}

In this section, we put forward some necessary preliminary lemmas, which are very important to the proofs of our main results.

For $t\leq\alpha(G)$, let $\delta_t(G)$ be the minimum order of the neighborhood of an independent set of order $t$ in a graph $G$. Namely,
$\delta_t(G)=\min\{|N_G(I)|: I \ \mbox{is an independent set of order} \ t\}$. Erbes, Molla, Mousley and Santana\cite{EMMS} proved the following
two results.

\medskip

\noindent{\textbf{Lemma 2.1}} (Erbes, Molla, Mousley and Santana \cite{EMMS}). Let $d\geq3$ be an integer, and let $G$ be a connected graph. Then $i(G-S)<\frac{2}{d-2}|S|$ for all nonempty $S\subseteq V(G)$ if and only if $\delta_t(G)>\frac{t(d-2)}{2}$ for all $t$ satisfying $1\leq t\leq\alpha(G)$.

\medskip

\noindent{\textbf{Lemma 2.2}} (Erbes, Molla, Mousley and Santana \cite{EMMS}). Let $d$ be an integer with $d\geq4$, and let $G$ be a connected graph
of order $n$ with $n>d$ and $\alpha(G)\leq5$. If
$$
\delta_{2t}(G)>t(d-2)
$$
for all $t$ satisfying $1\leq t\leq\frac{\alpha(G)}{2}$, then $G$ has a spanning tree with leaf distance at least $d$.

\medskip

Ma and Liu \cite{ML} showed a characterization for a graph to be fractional $k$-extendable.

\medskip

\noindent{\textbf{Lemma 2.3}} (\cite{ML}). Let $k\geq1$ be an integer, and let $G$ be a graph with a $k$-matching. Then $G$ is fractional
$k$-extendable if and only if
$$
i(G-S)\leq|S|-2k
$$
holds for any $S\subseteq V(G)$ such that $G[S]$ contains a $k$-matching.

\medskip

\noindent{\textbf{Lemma 2.4}} (\cite{B}). Let $G$ be a connected graph, and let $H$ be a proper subgraph of $G$. Then $\rho(G)>\rho(H)$.

\medskip

\noindent{\textbf{Lemma 2.5}} (Hong \cite{Ha}). Let $G$ be a graph with $n$ vertices. Then
$$
\rho(G)\leq\sqrt{2e(G)-n+1},
$$
where the equality holds if and only if $G$ is a star or a complete graph.

\medskip

Let $M$ be a real symmetric matrix whose rows and columns are indexed by $V=\{1,2,\cdots,n\}$.
Suppose that $M$ can be written as
\begin{align*}
M=\left(
  \begin{array}{ccc}
    M_{11} & \cdots & M_{1s}\\
    \vdots & \ddots & \vdots\\
    M_{s1} & \cdots & M_{ss}\\
  \end{array}
\right)
\end{align*}
in terms of partition $\pi: V=V_1\cup V_2\cup\cdots\cup V_s$ , wherein $M_{ij}$ is the submatrix (block) of $M$ obtained by rows in $V_i$ and
columns in $V_j$. The average row sum of $M_{ij}$ is denoted by $q_{ij}$. Then matrix $M_{\pi}=(q_{ij})$ is said to be the quotient matrix of $M$.
If the row sum of every block $M_{ij}$ is a constant, then the partition is equitable.

\medskip

\noindent{\textbf{Lemma 2.6}} (\cite{YYSX}). Let $M$ be a real matrix with an equitable partition $\pi$, and let $M_{\pi}$ be the corresponding
quotient matrix. Then every eigenvalue of $M_{\pi}$ is an eigenvalue of $M$. Furthermore, if $M$ is nonnegative, then the largest eigenvalues of
$M$ and $M_{\pi}$ are equal.

\medskip

The subsequent lemma is the well-known Cauchy Interlacing Theorem.

\medskip

\noindent{\textbf{Lemma 2.7}} (Haemers \cite{Hi}). Let $M$ be a Hermitian matrix of order $s$, and let $N$ be a principal submatrix of $M$ with
order $t$. If $\lambda_1\geq\lambda_2\geq\cdots\geq\lambda_s$ are the eigenvalues of $M$ and $\mu_1\geq\mu_2\geq\cdots\geq\mu_t$ are the eigenvalues
of $N$, then $\lambda_i\geq\mu_i\geq\lambda_{s-t+i}$ for $1\leq i\leq t$.

\section{The proof of Theorem 1.1}

In order to verify Theorem 1.1, we first prove the following lemma.

\medskip

\noindent{\textbf{Lemma 3.1.}} Let $d$ is an integer with $d\geq3$, and let $G$ be a connected graph of order $n$ with $n\geq d^{2}$. If
$$
\rho(G)\geq\rho(K_{\lceil\frac{d}{2}\rceil-1}\vee(K_{n-\lceil\frac{d}{2}\rceil}\cup K_1)),
$$
then $\delta_t(G)>\frac{t(d-2)}{2}$ for all $t$ satisfying $1\leq t\leq\alpha(G)$, unless $G=K_{\lceil\frac{d}{2}\rceil-1}\vee(K_{n-\lceil\frac{d}{2}\rceil}\cup K_1)$.

\medskip

\medskip

\noindent{\it Proof.} Suppose that $\delta_t(G)\leq\frac{t(d-2)}{2}$ for some $t$ satisfying $1\leq t\leq\alpha(G)$. According to Lemma 2.1,
we conclude
$$
\frac{(d-2)\cdot i(G-S)}{2}\geq|S|
$$
for some nonempty $S\subseteq V(G)$. By the integrity of $|S|$, we see
$$
\left\lceil\frac{(d-2)\cdot i(G-S)}{2}\right\rceil\geq|S|
$$
for some nonempty $S\subseteq V(G)$. Let $|S|=s$ and $i(G-S)=q$. Then $G$ is a spanning subgraph of
$G_1=K_{\lceil\frac{q(d-2)}{2}\rceil}\vee(K_{n_1}\cup qK_1)$, where $n_1=n-\lceil\frac{q(d-2)}{2}\rceil-q$. In view of Lemma 2.4, we obtain
\begin{align}\label{eq:3.1}
\rho(G)\leq\rho(K_{\lceil\frac{q(d-2)}{2}\rceil}\vee(K_{n-\lceil\frac{q(d-2)}{2}\rceil-q}\cup qK_1)),
\end{align}
where the equality holds if and only if $G=K_{\lceil\frac{q(d-2)}{2}\rceil}\vee(K_{n-\lceil\frac{q(d-2)}{2}\rceil-q}\cup qK_1)$.

If $q=1$, then $G_1=K_{\lceil\frac{d}{2}\rceil-1}\vee(K_{n-\lceil\frac{d}{2}\rceil}\cup K_1)$. Using \eqref{eq:3.1}, we get
$$
\rho(G)\leq\rho(K_{\lceil\frac{d}{2}\rceil-1}\vee(K_{n-\lceil\frac{d}{2}\rceil}\cup K_1)),
$$
with equality if and only if $G=K_{\lceil\frac{d}{2}\rceil-1}\vee(K_{n-\lceil\frac{d}{2}\rceil}\cup K_1)$, a contradiction. In what follows,
we consider $q\geq2$.

Recall that $G_1=K_{\lceil\frac{q(d-2)}{2}\rceil}\vee(K_{n-\lceil\frac{q(d-2)}{2}\rceil-q}\cup qK_1)$. By a direct computation, we have
\begin{align}\label{eq:3.2}
e(G_1)=&\binom{n-q}{2}+q\left\lceil\frac{q(d-2)}{2}\right\rceil\nonumber\\
=&\frac{(n-q)(n-q-1)}{2}+q\left\lceil\frac{q(d-2)}{2}\right\rceil.
\end{align}

It follows from \eqref{eq:3.2} and Lemma 2.5 that
\begin{align}\label{eq:3.3}
\rho(G_1)\leq&\sqrt{2e(G_1)-n+1}\nonumber\\
=&\sqrt{(n-q)(n-q-1)+2q\left\lceil\frac{q(d-2)}{2}\right\rceil-n+1}\nonumber\\
\leq&\sqrt{(n-q)(n-q-1)+2q\cdot\frac{q(d-2)+1}{2}-n+1}\nonumber\\
=&\sqrt{(d-1)q^{2}-(2n-2)q+n^{2}-2n+1}.
\end{align}
Let $\psi_1(q)=(d-1)q^{2}-(2n-2)q+n^{2}-2n+1$. Notice that $n\geq\lceil\frac{q(d-2)}{2}\rceil+q\geq\frac{q(d-2)}{2}+q=\frac{qd}{2}$. Then we
obtain $2\leq q\leq\frac{2n}{d}$. By a direct calculation, we have
$$
\psi_1\left(\frac{2n}{d}\right)-\psi_1(2)=-\frac{4}{d^{2}}(n-d)(n-d^{2})\leq0
$$
by $n\geq d^{2}$. Thus, we see that $\psi_1(q)$ attains its maximum value at $q=2$ for $2\leq q\leq\frac{2n}{d}$. Together with \eqref{eq:3.3},
$n\geq d^{2}$ and $d\geq3$, we get
\begin{align}\label{eq:3.4}
\rho(G_1)\leq&\sqrt{\psi_1(2)}\nonumber\\
=&\sqrt{4(d-1)-2(2n-2)+n^{2}-2n+1}\nonumber\\
=&\sqrt{(n-2)^{2}-2n+4d-3}\nonumber\\
\leq&\sqrt{(n-2)^{2}-2d^{2}+4d-3}\nonumber\\
=&\sqrt{(n-2)^{2}-2(d-1)^{2}-1}\nonumber\\
<&n-2.
\end{align}

Since $K_{n-1}$ is a proper subgraph of $K_{\lceil\frac{d}{2}\rceil-1}\vee(K_{n-\lceil\frac{d}{2}\rceil}\cup K_1)$, it follows from Lemma 2.4
that
\begin{align}\label{eq:3.5}
\rho(K_{\lceil\frac{d}{2}\rceil-1}\vee(K_{n-\lceil\frac{d}{2}\rceil}\cup K_1))>\rho(K_{n-1})=n-2.
\end{align}

Using \eqref{eq:3.1}, \eqref{eq:3.4} and \eqref{eq:3.5}, we conclude
$$
\rho(G)\leq\rho(G_1)<n-2<\rho(K_{\lceil\frac{d}{2}\rceil-1}\vee(K_{n-\lceil\frac{d}{2}\rceil}\cup K_1)),
$$
which contradicts $\rho(G)\geq\rho(K_{\lceil\frac{d}{2}\rceil-1}\vee(K_{n-\lceil\frac{d}{2}\rceil}\cup K_1))$. This completes the proof of
Lemma 3.1. \hfill $\Box$

\medskip

Next, we prove Theorem 1.1.

\medskip

\noindent{\it Proof of Theorem 1.1.} According to Lemma 3.1, we see
$$
\delta_k(G)>\frac{k(d-2)}{2}
$$
for all even $k$ satisfying $2\leq k\leq\alpha(G)$, unless $G=K_{\lceil\frac{d}{2}\rceil-1}\vee(K_{n-\lceil\frac{d}{2}\rceil}\cup K_1)$. Let
$k=2t$. Then we have
$$
\delta_{2t}(G)>t(d-2)
$$
for all $t$ satisfying $1\leq t\leq\frac{\alpha(G)}{2}$, unless $G=K_{\lceil\frac{d}{2}\rceil-1}\vee(K_{n-\lceil\frac{d}{2}\rceil}\cup K_1)$.
Combining this with $n\geq d^{2}\geq16$, $\alpha(G)\leq5$ and Lemma 2.2, we see that $G$ has a spanning tree with leaf distance at least $d$,
unless $G=K_{\lceil\frac{d}{2}\rceil-1}\vee(K_{n-\lceil\frac{d}{2}\rceil}\cup K_1)$. Theorem 1.1 is proved. \hfill $\Box$

\section{The proof of Theorem 1.2}

In this section, we prove Theorem 1.2.

\medskip

\noindent{\it Proof of Theorem 1.2.} Suppose, to the contrary, that $G$ is not fractional $k$-extendable. According to Lemma 2.3, there exists
some nonempty subset $S$ of $V(G)$ such that $|S|\geq2k$ and $i(G-S)\geq|S|-2k+1$. Then $G$ is a spanning subgraph of
$G_1=K_s\vee(K_{n_1}\cup(s-2k+1)K_1)$, where $|S|=s\geq2k$ and $n_1=n-2s+2k-1\geq0$. Using Lemma 2.4, we conclude
\begin{align}\label{eq:4.1}
\rho(G)\leq\rho(G_1),
\end{align}
with equality if and only if $G=G_1$. Notice that $G$ has the minimum degree $\delta$. Thus, we have $\delta(G_1)=s\geq\delta(G)=\delta$. Then
we proceed by the following two cases.

\medskip

\noindent{\bf Case 1.} $\delta\leq2k$.

Obviously, $s\geq2k\geq\delta$. Let $G_2=K_{2k}\vee(K_{n-2k-1}\cup K_1)$. We are to prove that $\rho(G_1)\leq\rho(G_2)$ with equality if and only
if $G_1=G_2$.

It is obvious that $G_1=G_2$ if $s=2k$, and so $\rho(G_1)=\rho(G_2)$. Next, we are to consider $s\geq2k+1$.

In terms of the partition $V(G_1)=V(K_s)\cup V(K_{n-2s+2k-1})\cup V((s-2k+1)K_1)$, the quotient matrix of $A(G_1)$ is equal to
\begin{align*}
B_1=\left(
  \begin{array}{ccc}
    s-1 & n-2s+2k-1 & s-2k+1\\
    s & n-2s+2k-2 & 0\\
    s & 0 & 0\\
  \end{array}
\right).
\end{align*}
Then the characteristic polynomial of the matrix $B_1$ is
\begin{align*}
\varphi_{B_1}(x)=x^{3}+(s-2k+3-n)x^{2}+(2ks-s^{2}-2k+2-n)x+s(s-2k+1)(n-2s+2k-2).
\end{align*}
Since the partition $V(G_1)=V(K_s)\cup V(K_{n-2s+2k-1})\cup V((s-2k+1)K_1)$ is equitable, it follows from Lemma 2.6 that $\rho(G_1)$ is the largest
root of $\varphi_{B_1}(x)=0$. Namely, $\varphi_{B_1}(\rho(G_1))=0$. Let $\gamma_1=\rho(G_1)\geq\gamma_2\geq\gamma_3$ be the three roots of $\varphi_{B_1}(x)=0$
and $Q_1=\mbox{diag}(s,n-2s+2k-1,s-2k+1)$. One checks that
\begin{align*}
Q_1^{\frac{1}{2}}B_1Q_1^{-\frac{1}{2}}=\left(
  \begin{array}{ccc}
    s-1 & s^{\frac{1}{2}}(n-2s+2k-1)^{\frac{1}{2}} & s^{\frac{1}{2}}(s-2k+1)^{\frac{1}{2}}\\
    s^{\frac{1}{2}}(n-2s+2k-1)^{\frac{1}{2}} & n-2s+2k-2 & 0\\
    s^{\frac{1}{2}}(s-2k+1)^{\frac{1}{2}} & 0 & 0\\
  \end{array}
\right)
\end{align*}
is symmetric, and also contains
\begin{align*}
\left(
  \begin{array}{ccc}
    n-2s+2k-2 & 0\\
    0 & 0\\
  \end{array}
\right)
\end{align*}
as its submatrix. Since $Q_1^{\frac{1}{2}}B_1Q_1^{-\frac{1}{2}}$ and $B_1$ have the same eigenvalues, by Lemma 2.7, we conclude
\begin{align}\label{eq:4.2}
\gamma_2\leq n-2s+2k-2<n-2.
\end{align}

Recall that $G_2=K_{2k}\vee(K_{n-2k-1}\cup K_1)$. Then the quotient matrix of $A(G_2)$ by the partition $V(G_2)=V(K_{2k})\cup V(K_{n-2k-1})\cup V(K_1)$
is equal to
\begin{align*}
B_2=\left(
  \begin{array}{ccc}
    2k-1 & n-2k-1 & 1\\
    2k & n-2k-2 & 0\\
    2k & 0 & 0\\
  \end{array}
\right),
\end{align*}
whose characteristic polynomial is
$$
\varphi_{B_2}(x)=x^{3}+(3-n)x^{2}+(2-2k-n)x+2k(n-2k-2).
$$
By virtue of Lemma 2.6, the largest root, say $\rho_2$, of $\varphi_{B_2}(x)=0$ is equal to $\rho(G_2)$.

Note that $K_{n-1}$ is a proper subgraph of $G_2=K_{2k}\vee(K_{n-2k-1}\cup K_1)$, it follows from \eqref{eq:4.2} and Lemma 2.4 that
\begin{align}\label{eq:4.3}
\rho_2=\rho(G_2)>\rho(K_{n-1})=n-2\geq\gamma_2.
\end{align}
Next, we prove $\varphi_{B_1}(\rho_2)=\varphi_{B_1}(\rho_2)-\varphi_{B_2}(\rho_2)>0$. By a direct calculation, we get
\begin{align}\label{eq:4.4}
\varphi_{B_1}(\rho_2)=\varphi_{B_1}(\rho_2)-\varphi_{B_2}(\rho_2)=(s-2k)f(\rho_2),
\end{align}
where $f(\rho_2)=\rho_2^{2}-s\rho_2+(s+1)n-2s^{2}+2ks-4s-2k-2$.

\noindent{\bf Claim 1.} $f(\rho_2)>0$ for $\rho_2>n-2$.

\noindent{\it Proof.} Firstly, we consider $n=2s-2k+1$. Together with $n\geq2k+9$, we deduce $s\geq2k+4$. Combining this with \eqref{eq:4.3}, we get
$$
\frac{s}{2}<2s-2k-1=n-2<\rho_2,
$$
and so
\begin{align*}
f(\rho_2)>&f(n-2)\\
=&n^{2}-3n-2s^{2}+2ks-2s-2k+2\\
=&(2s-2k+1)^{2}-3(2s-2k+1)-2s^{2}+2ks-2s-2k+2\\
=&2s^{2}-(6k+4)s+4k^{2}\\
\geq&2(2k+4)^{2}-(6k+4)(2k+4)+4k^{2}\\
=&16>0.
\end{align*}

Now we consider $n\geq2s-2k+2$. If $s\geq2k+2$, then it follows from \eqref{eq:4.3} that
$$
\frac{s}{2}<2s-2k\leq n-2<\rho_2,
$$
and so
\begin{align*}
f(\rho_2)>&f(n-2)\\
=&n^{2}-3n-2s^{2}+2ks-2s-2k+2\\
\geq&(2s-2k+2)^{2}-3(2s-2k+2)-2s^{2}+2ks-2s-2k+2\\
=&2s^{2}-6ks+4k^{2}-4k\\
\geq&2(2k+2)^{2}-6k(2k+2)+4k^{2}-4k\\
=&8>0.
\end{align*}
If $s=2k+1$, then we deduce
$$
\frac{s}{2}<2s-2k\leq n-2<\rho_2
$$
by \eqref{eq:4.3} and $n\geq2s-2k+2$. Combining this with $n\geq2k+9$, we conclude
\begin{align*}
f(\rho_2)>&f(n-2)\\
=&n^{2}-3n-2s^{2}+2ks-2s-2k+2\\
\geq&(2k+9)^{2}-3(2k+9)-2(2k+1)^{2}+2k(2k+1)-2(2k+1)-2k+2\\
=&8k+52>0.
\end{align*}
Claim 1 is proved. \hfill $\Box$

According to \eqref{eq:4.3}, \eqref{eq:4.4}, $s\geq2k+1$ and Claim 1, we obtain
$$
\varphi_{B_1}(\rho_2)=(s-2k)f(\rho_2)>0.
$$
As $\gamma_2\leq n-2<\rho(G_2)=\rho_2$ (see \eqref{eq:4.3}), we deduce
$$
\rho(G_1)=\gamma_1<\rho_2=\rho(G_2).
$$

From the above discussion, we have
\begin{align}\label{eq:4.5}
\rho(G_1)\leq\rho(G_2),
\end{align}
with equality if and only if $G_1=G_2$. Recall that $G_2=K_{2k}\vee(K_{n-2k-1}\cup K_1)$. It follows from \eqref{eq:4.1} and \eqref{eq:4.5} that
$$
\rho(G)\leq\rho(K_{2k}\vee(K_{n-2k-1}\cup K_1)),
$$
with equality if and only if $G=K_{2k}\vee(K_{n-2k-1}\cup K_1)$, a contradiction.

\medskip

\noindent{\bf Case 2.} $\delta\geq2k+1$.

Clearly, $s\geq\delta\geq2k+1$. Recall that $G_1=K_s\vee(K_{n-2s+2k-1}\cup(s-2k+1)K_1)$, the adjacency matrix $A(G_1)$ of $G_1$ has the quotient
matrix $B_1$, and $B_1$ has the characteristic polynomial $\varphi_{B_1}(x)$. Let $G_3=K_{\delta}\vee(K_{n-2\delta+2k-1}\cup(\delta-2k+1)K_1)$,
where $n\geq2\delta-2k+1$. We are to verify $\rho(G_1)\leq\rho(G_3)$ with equality if and only if $G_1=G_3$.

It is clear that $G_1=G_3$ if $s=\delta$, and so $\rho(G_1)=\rho(G_3)$. In what follows, we are to consider $s\geq\delta+1$.

For the graph $G_3$, its adjacency matrix $A(G_3)$ has the quotient matrix $B_3$ which is formed by replacing $s$ with $\delta$ in $B_1$, and $B_3$
has the characteristic polynomial $\varphi_{B_3}(x)$ which is derived by replacing $s$ with $\delta$ in $\varphi_{B_1}(x)$. Thus, we obtain
\begin{align*}
\varphi_{B_3}(x)=x^{3}+(\delta-2k+3-n)x^{2}+(2k\delta-\delta^{2}-2k+2-n)x+\delta(\delta-2k+1)(n-2\delta+2k-2).
\end{align*}
In terms of Lemma 2.6, the largest root, say $\rho_3$, of $\varphi_{B_3}(x)=0$ equals the spectral radius of $G_3$. That is, $\rho(G_3)=\rho_3$.

Notice that $\varphi_{B_3}(\rho_3)=0$. By plugging the value $\rho_3$ into $x$ of $\varphi_{B_1}(x)-\varphi_{B_3}(x)$, we have
\begin{align}\label{eq:4.6}
\varphi_{B_1}(\rho_3)=\varphi_{B_1}(\rho_3)-\varphi_{B_3}(\rho_3)=(s-\delta)g(\rho_3),
\end{align}
where $g(\rho_3)=\rho_3^{2}-(s+\delta-2k)\rho_3-2s^{2}+ns+6ks-2\delta s-4s-2kn+\delta n+n-2\delta^{2}+6k \delta-4\delta-4k^{2}+6k-2$. Since
$K_{n-\delta+k-1}$ is a proper subgraph of $G_3$, it follows from Lemma 2.4 that
\begin{align}\label{eq:4.7}
\rho_3=\rho(G_3)>\rho(K_{n-\delta+k-1})=n-\delta+k-2.
\end{align}
From \eqref{eq:4.7}, $s\geq\delta+1$ and $n\geq2s-2k+1$, we deduce
\begin{align}\label{eq:4.8}
\frac{s+\delta-2k}{2}<n-\delta+k-2<\rho_3,
\end{align}
which leads to
\begin{align}\label{eq:4.9}
g(\rho_3)>&g(n-\delta+k-2)\nonumber\\
=&n^{2}-(2\delta-2k+3)n-2s^{2}+(5k-\delta-2)s+k\delta+2\delta-k^{2}-2k+2.
\end{align}

We first consider $s\geq\frac{5}{2}\delta+1$. In view of \eqref{eq:4.9} and $n\geq2s-2k+1$, we conclude
\begin{align}\label{eq:4.10}
g(\rho_3)>&n^{2}-(2\delta-2k+3)n-2s^{2}+(5k-\delta-2)s+k\delta+2\delta-k^{2}-2k+2\nonumber\\
\geq&(2s-2k+1)^{2}-(2\delta-2k+3)(2s-2k+1)\nonumber\\
&-2s^{2}+(5k-\delta-2)s+k\delta+2\delta-k^{2}-2k+2\nonumber\\
=&2s^{2}-(5\delta-k+4)s+5k\delta-k^{2}+2k.
\end{align}
Let $h(s)=2s^{2}-(5\delta-k+4)s+5k\delta-k^{2}+2k$. Note that $\delta\geq2k+1$ and
$$
\frac{5\delta-k+4}{4}<\frac{5}{2}\delta+1\leq s,
$$
which implies that
\begin{align*}
h(s)\geq&h\left(\frac{5}{2}\delta+1\right)\\
=&\left(\frac{15}{2}k-5\right)\delta-k^{2}+3k-2\\
\geq&\frac{5}{2}k\delta-k^{2}+3k-2\\
\geq&\frac{5}{2}k(2k+1)-k^{2}+3k-2\\
=&4k^{2}+\frac{11}{2}k-2\\
>&0.
\end{align*}
Combining this with \eqref{eq:4.6}, \eqref{eq:4.10} and $s\geq\frac{5}{2}\delta+1$, we infer
\begin{align}\label{eq:4.11}
\varphi_{B_1}(\rho_3)=(s-\delta)g(\rho_3)>(s-\delta)h(s)>0.
\end{align}

According to \eqref{eq:4.7}, $s\geq\frac{5}{2}\delta+1$ and $n\geq2s-2k+1$, we have
\begin{align}\label{eq:4.12}
\varphi_{B_1}'(\rho_3)=&(s-\delta)g'(\rho_3)\nonumber\\
=&(s-\delta)(2\rho_3-s-\delta+2k)\nonumber\\
>&(s-\delta)(2(n-\delta+k-2)-s-\delta+2k)\nonumber\\
=&(s-\delta)(2n-3\delta-s+4k-4)\nonumber\\
\geq&(s-\delta)(2(2s-2k+1)-3\delta-s+4k-4)\nonumber\\
=&(s-\delta)(3s-3\delta-2)\nonumber\\
>&0.
\end{align}
The inequalities \eqref{eq:4.11} and \eqref{eq:4.12} imply
$$
\rho(G_1)=\gamma_1<\rho_3=\rho(G_3).
$$

From the above discussion, we get
\begin{align}\label{eq:4.13}
\rho(G_1)\leq\rho(G_3),
\end{align}
with equality if and only if $G_1=G_3$. Recall that $G_3=K_{\delta}\vee(K_{n-2\delta+2k-1}\cup(\delta-2k+1)K_1)$. By virtue of \eqref{eq:4.1}
and \eqref{eq:4.13}, we deduce
$$
\rho(G)\leq\rho(K_{\delta}\vee(K_{n-2\delta+2k-1}\cup(\delta-2k+1)K_1)),
$$
with equality if and only if $G=K_{\delta}\vee(K_{n-2\delta+2k-1}\cup(\delta-2k+1)K_1)$, a contradiction.

In what follows, we consider $\delta+1\leq s<\frac{5}{2}\delta+1$. Notice that $\frac{5k-\delta-2}{4}<\delta+1\leq s<\frac{5}{2}\delta+1$.
According to \eqref{eq:4.9}, $\delta\geq2k+1$ and $n\geq5\delta+1$, we obtain
\begin{align}\label{eq:4.14}
g(\rho_3)>&n^{2}-(2\delta-2k+3)n-2s^{2}+(5k-\delta-2)s+k\delta+2\delta-k^{2}-2k+2\nonumber\\
>&n^{2}-(2\delta-2k+3)n-2\left(\frac{5}{2}\delta+1\right)^{2}+(5k-\delta-2)\left(\frac{5}{2}\delta+1\right)\nonumber\\
&+k\delta+2\delta-k^{2}-2k+2\nonumber\\
=&n^{2}-(2\delta-2k+3)n-15\delta^{2}+\frac{27}{2}k\delta-14\delta-k^{2}+3k-2\nonumber\\
\geq&(5\delta+1)^{2}-(2\delta-2k+3)(5\delta+1)-15\delta^{2}+\frac{27}{2}k\delta-14\delta-k^{2}+3k-2\nonumber\\
=&\left(\frac{47}{2}k-21\right)\delta-k^{2}+5k-4\nonumber\\
\geq&\left(\frac{47}{2}k-21\right)(2k+1)-k^{2}+5k-4\nonumber\\
=&46k^{2}-\frac{27}{2}k-25\nonumber\\
>&0.
\end{align}
According to \eqref{eq:4.6}, \eqref{eq:4.14} and $\delta+1\leq s<\frac{5}{2}\delta+1$, we have
\begin{align}\label{eq:4.15}
\varphi_{B_1}(\rho_3)=(s-\delta)g(\rho_3)>0.
\end{align}

Using \eqref{eq:4.7}, $\delta+1\leq s<\frac{5}{2}\delta+1$ and $n\geq5\delta+1$, we obtain
\begin{align}\label{eq:4.16}
\varphi_{B_1}'(\rho_3)=&(s-\delta)g'(\rho_3)\nonumber\\
=&(s-\delta)(2\rho_3-s-\delta+2k)\nonumber\\
>&(s-\delta)(2(n-\delta+k-2)-s-\delta+2k)\nonumber\\
=&(s-\delta)(2n-3\delta-s+4k-4)\nonumber\\
\geq&(s-\delta)(2(5\delta+1)-3\delta-s+4k-4)\nonumber\\
=&(s-\delta)(7\delta-s+4k-2)\nonumber\\
>&0.
\end{align}
The inequalities \eqref{eq:4.15} and \eqref{eq:4.16} yield
$$
\rho(G_1)=\gamma_1<\rho_3=\rho(G_3).
$$

From the above discussion, we conclude
\begin{align}\label{eq:4.17}
\rho(G_1)\leq\rho(G_3),
\end{align}
with equality if and only if $G_1=G_3$. Recall that $G_3=K_{\delta}\vee(K_{n-2\delta+2k-1}\cup(\delta-2k+1)K_1)$. It follows from \eqref{eq:4.1}
and \eqref{eq:4.17} that
$$
\rho(G)\leq\rho(K_{\delta}\vee(K_{n-2\delta+2k-1}\cup(\delta-2k+1)K_1)),
$$
with equality if and only if $G=K_{\delta}\vee(K_{n-2\delta+2k-1}\cup(\delta-2k+1)K_1)$, a contradiction. This completes the proof of Theorem 1.2. \hfill $\Box$

\medskip

\section*{Data availability statement}

My manuscript has no associated data.

\section*{Declaration of competing interest}

The authors declare that they have no conflicts of interest to this work.

\section*{Acknowledgments}


This work was supported by the Natural Science Foundation of Jiangsu Province (Grant No. BK20241949). Project ZR2023MA078 supported by Shandong
Provincial Natural Science Foundation.

\end{document}